\documentclass[11pt]{amsart}

\usepackage{amsmath}
\usepackage{amssymb}
\usepackage{amsfonts}
\usepackage{geometry}
\usepackage{url}
\usepackage{color}
\usepackage{blkarray}
\usepackage{comment}
\usepackage{blindtext,rotating}
\usepackage{setspace}
\usepackage{amsthm}
\usepackage[all,cmtip]{xy}
\usepackage{amscd}
\usepackage{fourier-orns}

\usepackage{amsthm}
\theoremstyle{plain}
\newtheorem*{theorem*}{Theorem}

\newcommand{\oh}{\mathcal{O}}
\newcommand{\el}{\mathrm{L}}
\newcommand{\m}{\mathrm{M}}

\newcommand{\s}{\mathfrak{sl}_2}

\newcommand{\Hom}{\operatorname{Hom}}

\newcommand{\Res}{\operatorname{Res}}

\newcommand{\Irr}{\operatorname{Irr}}

\begin{document}

\title{The finite-dimensional representations of the rational Cherednik algebra of $E_8$ when $c=1/3$}

\author{Emily Norton}

\begin{abstract} We finish the classification of finite-dimensional irreducible representations of rational Cherednik algebras $H_c(W)$ with equal parameters for exceptional Coxeter groups by resolving the last open case: when $W=E_8$ and the denominator of $c$ is $3$.
\end{abstract}

\maketitle

The finite-dimensional representations in Category $\oh_c(W)$ of a rational Cherednik algebra $H_c(W)$ are exactly the cuspidals: $\el(\tau)\in\oh_c(W)$ is finite-dimensional if and only if $\Res^{\oh_c(W)}_{\oh_c(W')}\el(\tau)=0$ for all parabolic subgroups $W'\subset W$ \cite{BE}. In this sense the classification of the finite-dimensional irreducible representations forms the foundation of the representation theory of $H_c(W)$. When $W$ is an exceptional Coxeter group, work by Chmutova \cite{Ch}, Balagovic-Puranik \cite{BaP}, the author \cite{N},\cite{N2}, Griffeth-Gusenbauer-Juteau-Lanini \cite{GGJL}, and Losev-Shelley-Abrahamson \cite{LSa} has completed 
this classification for equal parameters $c$ except when $W=E_8$ and the denominator of $c$ is equal to $3$. 

For the rest of this note, $c=r/3$, $r\in\mathbb{N}$ and $\mathrm{gcd}(r,3)=1$, and we work in Category $\oh_c(E_8)$, which contains all finite-dimensional representations. Let $\m(\tau)\in\oh_c(E_8)$ be the Verma module of lowest weight $\tau\in\Irr E_8$, and let $\el(\tau)$ be its simple head. Let $V=\mathbb{C}^8$ be the reflection representation of $E_8$. According to Section 5.10 of \cite{GGJL}, if $\el(\tau)$ is finite-dimensional then: 
$$\tau\in\{\mathrm{triv},V,\phi_{28,8},\phi_{35,2},\phi_{50,8},\phi_{160,7},\phi_{175,12},\phi_{300,8},\phi_{840,13}\}\subset\Irr E_8.$$
Switching to the notation of \cite{GP} and \cite{GJ} for $\Irr E_8$, if $\el(\tau)$ is finite-dimensional then:
$$\tau\in\{1_x,8_z,28_x,35_x,50_x,160_z,175_x,300_x,840_z\}\subset\Irr E_8.$$
The recent paper of Losev and Shelley-Abrahamson \cite{LSa} counts that there are exactly eight finite-dimensional irreducibles. To finish the classification it suffices to identify for which $\tau$ in the list, $\el(\tau)$ is infinite-dimensional.

\begin{theorem*} Let $c=1/3$. Then $\el(50_x)$ is infinite-dimensional.
\end{theorem*}

\begin{proof}
When $W$ is a finite Coxeter group, there is a well-known $\s$-triple $(\mathbf{e},\mathbf{h},\mathbf{f})$ in $H_c(W)$, and any finite-dimensional $H_c(W)$-representation is also a finite-dimensional $\s$-representation \cite{GGOR}. The element $\mathbf{h}:=\sum_{i=1}^{\dim V} x_iy_i+y_ix_i$, where $\{x_i\}$ is an orthonormal basis for $V^*$ and $\{y_i\}$ is a dual basis for $V$, commutes with $W$, acts on the lowest weight $\tau$ of any irreducible $\el=\el(\tau)$ by a scalar $\mathbf{h}_c(\tau)$, and puts a $\mathbb{Z}$-grading on $\el$ \cite{GGOR}. If $\el$ is finite-dimensional then the equality $\dim_{\mathbb{C}}\el[-i]=\dim_{\mathbb{C}}\el[i]$ must hold for any $i\in\mathbb{Z}$, where $\el[i]$ denotes the graded piece of $\el$ in degree $i$, by representation theory of $\s$. We will show that $\dim_{\mathbb{C}}\el(50_x)[-1]<\dim_{\mathbb{C}}\el(50_x)[1]$, which means that $\el(50_x)$ cannot be finite-dimensional.

When $c=1/3$, $\mathbf{h}_c(\tau)=(\dim V)/2-c\sum_{\hbox{reflections }s\in E_8}s|_{\tau}=4-(1/3)120\tau(s)/\dim\tau$. Table C.6 of \cite{GP} gives the values of $\tau(s)$, from which $\mathbf{h}_c(50_x)=-12$. Next, the decomposition matrix of the Hecke algebra at a root of unity for $e=3$ and $W=E_8$ given in Table 7.15 of \cite{GJ} is a submatrix of the decomposition matrix for $\oh_{1/3}(E_8)$, but one must tensor the labels $\tau$ for rows and columns by the sign representation. Write $\tau':=\tau\otimes\mathrm{sign}$. Look at the column for $50_x$ in Table 7.15 of \cite{GJ}, which corresponds to the column for $50_x'$ in $\oh_{1/3}(E_8)$, and observe that every entry in the column $50_x$ of Table 7.15 is zero (excepting row $50_x$ itself of course) until the row labeled $700_{xx}$, where there is a $1$. Applying Lemma 3.6 from \cite{N}, this means that $\dim\Hom(\m(50_x'),\m(\tau'))=0$ for all $\tau\in\Irr E_8$, $\mathbf{h}_c(50_x)<\mathbf{h}_c(\tau)<\mathbf{h}_c(700_{xx})$, but $\dim\Hom(\m(50_x'),\m(700_{xx}'))=1$. On the other hand, $\mathbf{h}_c(700_{xx})=0$ and for all other $\sigma$ such that $\mathbf{h}_c(\sigma)=0$, it follows from Table 7.15 of \cite{GJ} and Lemma 3.6 of \cite{N} that $\dim\Hom(\m(50_x'),\m(\sigma'))=0$. Moreover, there is no $\tau\in\Irr E_8$ with $\mathbf{h}_c(\tau)=1$. By Lemma 3.5 of \cite{N}, it follows that for all $700_{xx}\neq\tau\in\Irr E_8$ satisfying $\mathbf{h}_c(50_x)<\mathbf{h}_c(\tau)< 2$, $\dim\Hom(\m(\tau),\m(50_x))=0$, while $\dim\Hom(\m(700_{xx}),\m(50_x))=1$. 

The expression in the Grothendieck group of $\oh_{1/3}(E_8)$ for $\el(50_x)$ in the basis of Vermas $\m(\tau)$ 
is therefore: $\el(50_x)=\m(50_x)-\m(700_{xx})+\sum_{\mathbf{h}_c(\tau)\geq2}a_{50_x,\tau}\m(\tau)$. The graded character of $\el(50_x)$ truncated after degree $1$ is then $\sum_{k=0}^{13}{7+k\choose 7}50t^{-12+k} -\sum_{j=0}^1 {7+j\choose 7}700t^j$, and so the graded dimensions of $\el(50_x)$ in degrees $-1$ and $1$ (the coefficients of $t^{-1}$ and $t$ in the graded character) are:
\begin{align*}
\dim\el(50_x)[-1]&={18\choose 7}50=1,591,200\\
\dim\el(50_x)[1]&={20\choose 7}50-{8\choose 7}700=3,870,400
\end{align*}
Since $\dim\el(50_x)[-1]<\dim\el(50_x)[1]$, $\el(50_x)$ must be infinite-dimensional.
\end{proof}
It follows that the finite-dimensional irreducible representations of $H_{1/3}(E_8)$ are $\el(\tau)$ where $\tau\in\{1_x,8_z,28_x,35_x,160_z,175_x,300_x,840_z\}\subset\Irr E_8$. By category equivalences, the cases of $c=r/3$ for $r\neq1$ reduce to the case $c=1/3$ \cite{GGOR}. This completes the classification. 

\end{document}